\documentclass[12pt,a4paper]{article}
\usepackage{mathrsfs}
\usepackage{amssymb}
\usepackage{amsmath}
\usepackage{amsfonts}
\usepackage{amstext}
\usepackage{amsthm}
\usepackage{graphicx}
\usepackage[top=1in, bottom=1in,left=0.75in,right=0.75in]{geometry}

\def \qed {\hfill \vrule height6pt width 6pt depth 0pt}

\begin{document}

\title{Unilateral global interval bifurcation theorem for $p$-Laplacian and its applications
\thanks{Research supported by the NSFC (No. 11061030).}}
\author{{\small  Guowei Dai\thanks{Corresponding author. Tel: +86 931
7971124.\newline
\text{\quad\,\,\, E-mail address}: daiguowei@nwnu.edu.cn (G. Dai).
}
} \\
{\small Department of Mathematics, Northwest Normal
University, Lanzhou, 730070, PR China}\\
}
\date{}
\maketitle

\begin{abstract}
In this paper, we establish a unilateral global bifurcation result from interval for a class of $p$-Laplacian problems.
By applying the above result, we study the spectrum of a class of half-quasilinear problems.
Moreover, we also investigate the existence of nodal solutions for a class of half-quasilinear eigenvalue problems.
\\ \\
\textbf{Keywords}: Unilateral interval bifurcation; Half-quasilinear problems; Nodal solutions; $p$-Laplacian
\\ \\
\textbf{MSC(2000)}: 34C10; 34C23
\end{abstract}\textbf{\ }

\numberwithin{equation}{section}

\numberwithin{equation}{section}

\section{Introduction}
\bigskip

\quad\, Let $E$ be a real Banach space with the norm $\Vert \cdot\Vert$. Consider the operator equation
\begin{equation}\label{oe}
u=\lambda Lu+H(\lambda,u),
\end{equation}
where $L$ is a compact linear operator and $H:\mathbb{R}\times E\rightarrow E$ is compact
with $H(\lambda,u)=o(\Vert u\Vert)$ at $u=0$ uniformly on bounded $\lambda$ intervals. In [\ref{K}], Krasnosel'skii has shown that
all characteristic values of $L$ which are of odd multiplicity are bifurcation points.
We refer to the books [\ref{L1}, Theorem 6.2.1] and [\ref{LM2}, Theorem 12.1.4] and their references for the very latest refinements of the classical result by Krsnosel'skii; they are
valid in much more general contexts.
Rabinowitz [\ref{R2}] has extended
this result by showing that bifurcation has global consequences. More precisely, if $\mu$ is of odd multiplicity and
\begin{equation}
\mathscr{S}:=\overline{\left\{(\lambda,u)\big|(\lambda,u)\,\text{\,satisfies (\ref{oe}) and $u\not\equiv 0$\,\,}\right\}}^{\mathbb{R}\times E},\nonumber
\end{equation}
then $\mathscr{S}$ possesses a component which contains $(\mu,0)$ and either unbounded or meets another characteristic value of $L$.
We note that the nonlinear term $\lambda Lu+H(\lambda,u)$ is Fr\'{e}chet differentiable at the origin.

In the celebrated work [\ref{R2}], Rabinowitz also established unilateral global bifurcation theory.
However, as pointed out by Dancer [\ref{D1}, \ref{D2}] and L\'{o}pez-G\'{o}mez [\ref{L1}], the proofs of these
theorems contain gaps. Fortunately, L\'{o}pez-G\'{o}mez gave a corrected version of the unilateral global bifurcation
theorem for (\ref{oe}) [\ref{L1}, Theorem 6.4.3]. This is the first available correct unilateral theorem. Later, working on the theory of L\'{o}pez-G\'{o}mez [\ref{L1}] Dancer got
another unilateral theorem in [\ref{D2}] which has been extended to the one-dimensional $p$-Laplacian problem by Dai and Ma [\ref{DM}].

In [\ref{Be}], Berestycki considered
a class of problems involving nondifferentiable nonlinearity. More precisely, he considered the nonlinear
Sturm-Liouville problem
\begin{equation}\label{SL1}
\mathscr{L}u=\lambda au+f\left(x,u,u',\lambda\right)+g\left(x,u,u',\lambda\right),\,\,x\in(0,\pi)
\end{equation}
together with the following separated boundary conditions
\begin{equation}
b_0u(0)+c_0u'(0)=0,\nonumber
\end{equation}
\begin{equation}
b_1u(\pi)+c_1u'(\pi)=0,\nonumber
\end{equation}
where $b_i$, $c_i$ are real numbers such that $\vert b_i\vert+\vert c_i\vert\neq 0$, $i=0,1$;
$\mathscr{L}$ is a regular Sturm-Liouville operator and $a$ is a continuous positive function defined on $[0,\pi]$. It is assumed that $f$ and $g$ are continuous, with $\vert f\vert\leq M\vert u\vert$ in a neighborhood of $u=u'=0$ for $\forall x\in[0,\pi]$ and $\forall \lambda\in \mathbb{R}$ and $g=o\left(\vert u\vert+\vert u'\vert\right)$ near $(u,u')=(0,0)$ uniformly in $x\in[0,\pi]$ and $\lambda\in \Lambda$ for every bounded interval $\Lambda$. If $\lambda_1<\lambda_2<\cdots<\lambda_n<\cdots$ are the eigenvalues of $\mathscr{L}$, Berestycki shows that at least one continuum of solutions bifurcates from each interval $\left[\lambda_k-M/a_0,\lambda_k+M/a_0\right]$ where $a_0=\min_{x\in[0,\pi]}a(x)$. Furthermore, at least one of the continua bifurcating from $\left[\lambda_k-M/a_0,\lambda_k+M/a_0\right]$ has the nodal property, i.e., there exists a continuum $\mathfrak{D}$ emanating from $\left[\lambda_k-M/a_0,\lambda_k+M/a_0\right]$ such that for any $(\lambda,u)\in \mathfrak{D}$ then $u$ has fixed number of simple zeros. The above result has been improved partially by Schmitt and Smith [\ref{SS}] by applying a set-valued version Rabinowitz's global bifurcation theorem. Clearly, Eq. (\ref{SL1}) does not in general have a linearized about $u=0$ because of the presence of the term $f$.

Of course, the natural question is whether the results of [\ref{Be}] can be extended to the case that $\mathscr{L}$ is the quasilinear $p$-Laplacian operator. Meanwhile, another question is whether the interval version of bifurcation results of [\ref{DM}] exist. In this paper, we shall provide a positive answer to these questions.

For simplicity we shall restrict ourselves to the 0-Dirichlet boundary value problems, but the methods
used in this paper are also suitable for the separated boundary value problems like of [\ref{Be}]. Concretely, we shall study the following $p$-Laplacian problem
\begin{equation}\label{dF}
\left\{
\begin{array}{l}
\left(\left\vert u'\right\vert^{p-2}u'\right)'+\lambda a(x)\varphi_p(u)+F(x,u,\lambda)=0, \,\, x\in(0,1),\\
u(0)=u(1)=0,
\end{array}
\right.
\end{equation}
where $1<p<+\infty$, $\varphi_p(u)=\vert u\vert^{p-2}u$, $\lambda$ is a parameter, $a(x)\in C[0,1]$ is positive and $F:(0,1)\times \mathbb{R}^2\rightarrow \mathbb{R}$ is a continuous function. Moreover, the nonlinear term $F$ has the form $F=f+g$, where $f$ and $g$ satisfy the conditions:\\

(C1) $\left\vert f(x,s,\lambda)\right\vert\leq M\left\vert \varphi_p(s)\right\vert$ for all $x\in(0,1)$, $0<\vert s\vert\leq1$ and $\lambda\in \mathbb{R}$, where $M$ is a positive constant;

(C2) $g(x,s,\lambda)=o\left(\vert s\vert^{p-1}\right)$ near $s=0$, uniformly in $x\in(0,1)$ and on bounded $\lambda$ intervals.\\

Let $E:=\left\{u\in C^1[0,1]\big|u(0)=u(1)=0\right\}$ with the norm
\begin{equation}
\Vert u\Vert=\max_{x\in[0,1]}\vert u(x)\vert+\max_{x\in[0,1]}\left\vert u'(x)\right\vert.\nonumber
\end{equation}
Let $S_k^+$ denote the set of functions
in $E$ which have exactly $k-1$ interior nodal (i.e., non-degenerate) zeros in (0,1) and are positive near $x=0$, and set $S_k^-=-S_k^+$, and
$S_k =S_k^+\cup S_k^-$. It is clear that $S_k^+$ and $S_k^-$ are disjoint and open in $E$.
We also let $\Phi_k^{\pm}=\mathbb{R}\times S_k^{\pm}$ and $\Phi_k=\mathbb{R}\times S_k$ under the
product topology. We use $\mathscr{S}$ to denote the closure of the nontrivial solutions set of problem (\ref{dF}) in $\mathbb{R}\times E$,
and $\mathscr{S}_k^\pm$ to denote the subset of $\mathscr{S}$ with $u\in S_k^{\pm}$ and $\mathscr{S}_k=\mathscr{S}_k^+\cup \mathscr{S}_k^-$.
Finally, let $\lambda_k$ be the $k$th eigenvalue of the following eigenvalue problem
\begin{equation}\label{pe}
\left\{
\begin{array}{l}
-\left(\left\vert u'\right\vert^{p-2}u'\right)'=\lambda a(x)\varphi_p(u), \,\, x\in(0,1),\\
u(0)=u(1)=0.
\end{array}
\right.
\end{equation}
\indent Our main result for problem (\ref{dF}) is the following.
\\ \\
\textbf{Theorem 1.1.} \emph{Let $d=M/a_0$, where $a_0=\min_{x\in[0,1]}a(x)$, and let
$I_k=\left[\lambda_k-d,\lambda_k+d\right]$ for every $k\in \mathbb{N}$. The component $\mathscr{D}_k^+$ of $\mathscr{S}_k^+\cup\left(I_k\times \{0\}\right)$, containing $I_k\times \{0\}$ is unbounded and lies in $\Phi_k^+\cup\left(I_k\times \{0\}\right)$ and the component $\mathscr{D}_k^-$ of $\mathscr{S}_k^-\cup\left(I_k\times \{0\}\right)$, containing $I_k\times \{0\}$ is unbounded and lies in $\Phi_k^-\cup\left(I_k\times \{0\}\right)$.}\\
\\
\indent Note that the proofs of Lemma 1 and Theorem 1 of [\ref{Be}] strictly depend on the linear property of the operator $\mathscr{L}$. Thus, the methods used in [\ref{Be}] can not be used here to deal with the quasilinear problems (\ref{dF}). We use the generalized Picone identity to overcome the difficulty which is raised by quasilinear operator.
Moreover, we use the unilateral global
bifurcation theorem of [\ref{DM}] rather than the global bifurcation theorem of [\ref{R2}] which
is used by Berestycki in [\ref{Be}] to prove Theorem 1.1. Hence, Theorem 1.1 improves the corresponding result of [\ref{Be}, Theorem 1] even in the case of $p=2$.

On the basis of the unilateral global interval bifurcation result,
we establish the spectrum of the following half-quasilinear problem
\begin{equation}\label{hqle}
\left\{
\begin{array}{l}
-\left(\varphi_p\left(u'\right)\right)'=\lambda a(x)\varphi_p(u)+ \alpha \varphi_p\left(u^+\right)+ \beta \varphi_p\left(u^-\right), \,\, x\in
(0,1),\\
u(0)=u(1)=0,
\end{array}
\right.
\end{equation}
where $u^+=\max\{u,0\}$, $u^-=-\min\{u,0\}$, $\alpha$ and $\beta$ are two continuous functions defined on $[0,1]$.
More precisely, we shall use Theorem 1.1 to prove the following result.
\\ \\
\textbf{Theorem 1.2.} \emph{There exist two sequences of simple half-eigenvalues for problem (\ref{hqle}),
$\lambda_1^+<\lambda_2^+<\cdots<\lambda_k^+<\cdots$ and $\lambda_1^-<\lambda_2^-<\cdots<\lambda_k^-<\cdots$.
The corresponding half-linear solutions are in $\left\{\lambda_k^+\right\}\times S_k^+$ and $\left\{\lambda_k^-\right\}\times S_k^-$.
Furthermore, aside from these solutions and the trivial ones, there is no other solutions of problem (\ref{hqle}).}
\\ \\
\indent Furthermore, following the above eigenvalue theory, we shall investigate the existence of nodal
solutions for the following $p$-Laplacian problem
\begin{equation}\label{dn}
\left\{
\begin{array}{l}
\left(\left\vert u'\right\vert^{p-2}u'\right)'+ra(x)f(u)+\alpha(x)\varphi_p\left(u^+\right)+\beta(x)\varphi_p\left(u^-\right)=0, \,\, x\in(0,1),\\
u(0)=u(1)=0,
\end{array}
\right.
\end{equation}
where $f:\mathbb{R}\rightarrow \mathbb{R}$ is a continuous function, $r$ is a real parameter. Throughout this paper, we assume that $f$ satisfies the conditions:\\

(A1) $sf(s)>0$ for $s\neq 0$;

(A2) there exist $f_0, f_{\infty}\in(0,+\infty)$ such that
\begin{equation}
f_0=\lim\limits_{|s|\rightarrow 0}\frac{f(s)}{\varphi_p(s)},\, \, f_{\infty}=\lim\limits_{|s|\rightarrow +\infty}\frac{f(s)}{\varphi_p(s)}.\nonumber
\end{equation}
\indent The last main result of this paper is the following:\\ \\
\textbf{Theorem 1.3.} \emph{Assume that $f$ satisfies (A1) and (A2), and for some $k\in\mathbb{N}$, $\nu=+$ and $-$, either
$\lambda_k^\nu/f_{\infty}<r<\lambda_k^\nu/f_0$ or $\lambda_k^\nu/f_0<r<\lambda_k^\nu/f_{\infty}.$ Then problem (\ref{dn}) has a solution $u_k^\nu$ such that $u_k^\nu$ has exactly $k-1$ zeros in (0,1) and $\nu u_k^\nu$ is positive near 0.}
\\ \\
\indent In the case of $p=2$ and $\alpha=\beta\equiv0$, Ma and Thompson [\ref{MT1}] considered problem (\ref{dn}) with determining
interval of $r$ by the bifurcation theory of Rabinowitz [\ref{R1}, \ref{R2}], in which there exist nodal solutions for problem (\ref{dn}) under the assumptions of (A1) and (A2).
We note that the assumption $f_0\in(0,+\infty)$ implies that $f$ is Fr\'{e}chet differentiable at the origin, i.e., $f$ is linearizable at the origin. Moreover, the Fr\'{e}chet derivative of $f$ at the point $u=0$ in the direction $v$ is $f_0v$.
In the case of $p=2$, $\alpha=\beta\equiv0$ but $a$ changes its sign, Hess and Kato [\ref{HK}]
proved some well-known classical results which show that the principal eigenvalues of the weighted boundary value problem (\ref{pe}) are bifurcation
points to positive solutions.
The idea of using bifurcation methods to study the solvability of nonlinear boundary value
problems also has been applied to study various boundary value problems, for instance, see [\ref{H}, \ref{L0}, \ref{L2}].

For $p\neq2$ but $\alpha=\beta\equiv0$, Dai and Ma [\ref{DM}] have established the existence of nodal solutions for problem (\ref{dn}) with crossing nonlinearity which extends the results of [\ref{MT1}]. We also note that the assumption of $\left(ra(x)f(u)\right)/\varphi_p(u)$ crossing eigenvalues implies that $f$ is $p-$1-homogeneous linearizable at the origin and infinity, i.e., $f_0$, $f_\infty\in(0,+\infty)$.
We also note that, in high-dimensional case, there are also a lot of fundamental papers on the global bifurcation for $p$-Laplacian [\ref{DPM}, \ref{DH}, \ref{FR}, \ref{FMT}, \ref{GS}, \ref{GT}].

In the previously mentioned papers, the nonlinearities are Fr\'{e}chet differentiable or $p-$1-homogeneous linearizable at the origin or infinity. However, the nonlinear term of problem (\ref{dn}) is not necessary $p-$1-homogeneous linearizable at the origin and infinity because of the influence of the term $\alpha(x)\varphi_p\left(u^+\right)+\beta(x)\varphi_p\left(u^-\right)$. So the bifurcation theory of [\ref{DM}, \ref{D1},  \ref{L1}, \ref{R1}, \ref{R2}, \ref{R3}] cannot be applied directly to obtain our results. Luckily, using Theorem 1.1 and 1.2, we can obtain
some results of the existence of nodal solutions which extend the corresponding ones of [\ref{DM}, \ref{MT1}] in some sense.

The rest of this paper is arranged as follows. In Section 2, we give the proof of Theorem 1.1. In Section 3, we shall prove Theorem 1.2; as a byproduct, it is also shown that for a problem possessing jumping nonlinearities, these half-eigenvalues correspond to bifurcation points in a unilateral global sense. Theorem 1.3 is proved in the last Section; in this section, we also give a nonexistence result for problem (\ref{dn}).

\section{Unilateral global bifurcation from interval}

\bigskip

\quad\, Now, we consider the operator equation (\ref{oe}) again.
Rabinowitz's global bifurcation theorem [\ref{R2}, \ref{R4}] has shown that if the characteristic value $\mu$ of $L$ is of odd multiplicity,
then there exists a component $\mathscr{C}_\mu$ of $\mathscr{S}$ which contains $(\mu,0)$ and either unbounded or meets another characteristic value of $L$. Moreover, if $\mu$ is simple, Dancer [\ref{D1}] has shown that there are two distinct unbounded sub-continua $\mathscr{C}_\mu^+$ and $\mathscr{C}_\mu^-$ of the continuum $\mathscr{C}_\mu$ from ($\mu, 0)$, which satisfy either $\mathscr{C}_\mu^+$ and $\mathscr{C}_\mu^-$ are both unbounded or $\mathscr{C}_\mu^+\cap\mathscr{C}_\mu^-\neq\{(\mu,0)\}$. The result has been extended to the one-dimensional $p$-Laplacian problem by Dai and Ma [\ref{DM}]. More specifically, Dai and Ma [\ref{DM}] considered the following one-dimensional $p$-Laplacian problem
\begin{equation}\label{eg1}
\left\{
\begin{array}{l}
-\left(\varphi_p\left(u'\right)\right)'=\mu m(x)\varphi_p(u)+g(x,u,\mu), \, \, \text{a.e.}\,\, x\in
(0,1),\\
u(0)=u(1)=0,
\end{array}
\right.
\end{equation}
where $\varphi_p(s)=\vert s\vert^{p-2}s$, $1<p<+\infty$, $\mu$ is a positive parameter, $m(x)\geq 0$
and $m(x)\not\equiv 0$ for $x\in(0,1)$ is a continuous weight function,
$g:(0,1)\times \mathbb{R}\times\mathbb{R}\rightarrow\mathbb{R}$ satisfies the Carath\'{e}odory condition and
\begin{equation}\label{hg1}
\lim_{ s\rightarrow0}\frac{g(x,s;\mu)}{\vert s\vert^{p-1}}=0
\end{equation}
uniformly for a.e. $x\in(0,1)$ and $\mu$ on bounded sets. Let $\mu_k$ be the
$k$th eigenvalue of the corresponding linear problem of problem (\ref{eg1}).

They have shown that there are two distinct unbounded sub-continua $\mathscr{C}_k^+$ and $\mathscr{C}_k^-$ of the continuum $\mathscr{C}_k$ of problem (\ref{eg1}) emanating from ($\mu_k, 0)$, which satisfy:
\\ \\
\textbf{Lemma 2.1} (Theorem 3.2, [\ref{DM}]). \emph{Let $\nu\in\{+,-\}$. Then $\mathscr{C}_{k}^\nu$ is unbounded in $\mathbb{R}\times E$ and
\begin{equation}
\mathscr{C}_{k}^\nu\subset \left(\{(\mu_k,0)\}\cup\left(\mathbb{R}\times S_k^\nu\right)\right)\,\,\text{or}\,\, \mathscr{C}_{k}^\nu\subset \left(\{(\mu_k,0)\}\cup\left(\mathbb{R}\times S_k^{-\nu}\right)\right).\nonumber
\end{equation}}
\indent Next, we show that the existence and uniqueness theorem is valid for problem (\ref{dF}).\\ \\
\textbf{Lemma 2.2.} \emph{If $(\lambda, u)$ is a solution of problem (\ref{dF}) under the assumptions of (C1) and (C2)
and $u$ has a double zero, then $u \equiv 0$.}
\\ \\
\textbf{Proof.} Let $u$ be a solution of problem (\ref{dF}) and $x^*\in[0, 1]$ be a double zero.
We note that
\begin{equation}
u(x)=\int_x^{x^*}\varphi_p^{-1}\left(\int_s^{x^*}\left(-\lambda
a(\tau)\varphi_p(u(\tau))-f(\tau,u(\tau),\lambda)-g(\tau,u(\tau),\lambda)\right)\,d\tau\right)\,ds.\nonumber
\end{equation}
Firstly, we consider $x\in[0, x^*]$. Then we have that
\begin{eqnarray}
\vert u(x)\vert&\leq&\int_x^{x^*}\varphi_p^{-1}\left(\left\vert\int_s^{x^*}\left(-\lambda
a(\tau)\varphi_p(u(\tau))-f(\tau,u(\tau),\lambda)-g(\tau,u(\tau),\lambda)\right)\,d\tau\right\vert\right)\,ds\nonumber\\
&\leq&\varphi_p^{-1}
\left(\int_x^{x^*}
 \left\vert\lambda a(\tau)\varphi_p(u(\tau))+f(\tau,u(\tau),\lambda)+g(\tau,u(\tau),\lambda)\right\vert\,d\tau\right),\nonumber
\end{eqnarray}
furthermore, it follows that
\begin{eqnarray}
\varphi_p(\vert u(x)\vert)&\leq&\int_x^{x^*}
 \left\vert\lambda a(\tau)\varphi_p(u(\tau))+f(\tau,u(\tau),\lambda)+g(\tau,u(\tau),\lambda)\right\vert\,d\tau \nonumber\\
&\leq&\int_x^{x^*}
 \left\vert\lambda a(\tau)
+\frac{f(\tau,u(\tau),\lambda)}{\varphi_p(u(\tau))}+\frac{g(\tau,u(\tau),\lambda)}{\varphi_p(u(\tau))}\right\vert\varphi_p(u(\tau))\,d\tau \nonumber\\
&\leq&\int_x^{x^*}
 \left(\lambda a(\tau)
+\left\vert\frac{f(\tau,u(\tau),\lambda)}{\varphi_p(u(\tau))}\right\vert+\left\vert\frac{g(\tau,u(\tau),\lambda)}{\varphi_p(u(\tau))}\right\vert\right)\varphi_p(\vert u(\tau)\vert)\,d\tau.\nonumber
\end{eqnarray}
In view of (C2), for any $\varepsilon>0$, there exists a constant $1\geq\delta>0$ such that
\begin{equation}
\vert g(x,s,\lambda)\vert\leq \varepsilon\varphi_p( \vert s\vert)\nonumber
\end{equation}
uniformly with respect to $x\in(0,1)$ and fixed $\lambda$ when $\vert s\vert\in[0,\delta]$. Hence, we get that
\begin{eqnarray}
\varphi_p(\vert u(x)\vert)\leq\int_x^{x^*}
G(\tau,\lambda)\varphi_p(\vert u(\tau)\vert)\,d\tau,\nonumber
\end{eqnarray}
where
\begin{eqnarray}
G(\tau,\lambda)=\vert \lambda\vert\max_{x\in[0,1]}a(t)+M+\max_{\vert s\vert\in\left[\delta,\Vert u\Vert_\infty\right]}\left\vert\frac{f(\tau,s,\lambda)}{\varphi_p(s)}\right\vert+\max_{\vert s\vert\in\left[\delta,\Vert u\Vert_\infty\right]}\left\vert\frac{g(\tau,s,\lambda)}{\varphi_p(s)}\right\vert.\nonumber
\end{eqnarray}
By the Gronwall-Bellman inequality [\ref{Bre}, \ref{E}], we get $u \equiv 0$ on $[0, x^*]$.
Similarly, using a modification of Gronwall-Bellman inequality [\ref{ILL}, Lemma 2.2], we can get $u \equiv 0$ on $[x^*, 1]$.\qed\\ \\
\textbf{Remark 2.1.} By Lemma 2.2, we can see that if $(\lambda, u)$ is a nontrivial solution of problem (\ref{dF}) under the assumptions of (C1) and (C2), then $u\in \cup_{k=1}^\infty S_k$.\\
\\
\indent To prove Theorem 1.1, we introduce the following approximate problem
\begin{equation}\label{dnf}
\left\{
\begin{array}{l}
\left(\left\vert u'\right\vert^{p-2}u'\right)'+\lambda a(x)\varphi_p(u)+f\left(x,u\vert u\vert^\varepsilon,\lambda\right)+g(x,u,\lambda)=0, \,\, x\in(0,1),\\
u(0)=u(1)=0.
\end{array}
\right.
\end{equation}
The next lemma will play a key role in this paper which provides uniform a priori bounds for the solutions of problem (\ref{dnf}) near the
trivial solutions and will also ensure that $\left(\mathscr{S}_k^\nu\cap\left(\mathbb{R}\times\{0\}\right)\right)\subset\left(I_k\times\{0\}\right)$.
\\ \\
\textbf{Lemma 2.3.} \emph{Let $\epsilon_n$, $0\leq \epsilon_n\leq 1$, be a sequence converging to 0. If there exists a sequence $\left(\lambda_n,u_n\right)\in \mathbb{R}\times S_k^\nu$ such that $\left(\lambda_n,u_n\right)$ is a nontrivial solution of problem (\ref{dnf}) corresponding to $\epsilon=\epsilon_n$, and $\left(\lambda_n,u_n\right)$ converges to $(\lambda,0)$ in $\mathbb{R}\times E$, then $\lambda\in I_k$.}
\\ \\
\textbf{Proof.} Let $w_n=u_n/\left\Vert u_n\right\Vert$, then $w_n$ should be a solution of the problem
\begin{equation}\label{dnn}
\left\{
\begin{array}{l}
-\left(\varphi_p\left(w_n'\right)\right)'=\lambda_n a(x)\varphi_p(w_n)+\frac{f\left(x,u_n\vert u_n\vert^{\varepsilon_n},\lambda_n\right)}{\left\Vert u_n\right\Vert^{p-1}}+\frac{g\left(x,u_n,\lambda_n\right)}{\left\Vert u_n\right\Vert^{p-1}}, \,\, x\in(0,1),\\
w_n(0)=w_n(1)=0.
\end{array}
\right.
\end{equation}
Let
\begin{equation}
\widetilde{g}(x,u,\lambda)=\max_{0\leq \vert s\vert\leq u}\vert g(x,s,\lambda)\vert\,\,
\text{for all}\,\, x\in(0,1) \text{\,\,and\,\,} \lambda \text{\,\,on bounded sets},\nonumber
\end{equation}
then $\widetilde{g}$ is nondecreasing with respect to $u$ and
\begin{equation}\label{eg0+}
\lim_{ u\rightarrow 0^+}\frac{\widetilde{g}(x,u,\lambda)}{
u^{p-1}}=0.
\end{equation}
Further it follows from (\ref{eg0+}) that
\begin{eqnarray}\label{egn0}
\left\vert\frac{g(x,u,\lambda)}{\Vert u\Vert^{p-1}}\right\vert &\leq&\frac{
\widetilde{g}(x,\vert u\vert,\lambda)}{\Vert u\Vert^{p-1}}
\leq\frac{\widetilde{g}\left(x,\Vert u\Vert_\infty,\lambda\right)}{\Vert u\Vert^{p-1}}\nonumber\\
&\leq&\frac{ \widetilde{g}(x,\Vert u\Vert,\lambda)}{\Vert
u\Vert^{p-1}}\rightarrow0\,\,  \text{as}\,\, \Vert
u\Vert\rightarrow 0
\end{eqnarray}
uniformly in $x\in(0,1)$ and $\lambda$ on bounded sets. Clearly, (C1) implies that
\begin{eqnarray}\label{egn1}
\left\vert\frac{f\left(x,u_n\left\vert u_n\right\vert^{\epsilon_n},\lambda_n\right)}{\left\Vert u_n\right\Vert^{p-1}}\right\vert
&=&\left\vert\frac{f\left(x,u_n\left\vert u_n\right\vert^{\epsilon_n},\lambda_n\right)}{\varphi_p\left(u_n\left\vert u_n\right\vert^{\epsilon_n}\right)}\frac{\varphi_p\left(u_n\left\vert u_n\right\vert^{\epsilon_n}\right)}{\left\Vert u_n\right\Vert^{p-1}}\right\vert\nonumber\\
&\leq& M\left\Vert u_n\right\Vert^{(p-1)\epsilon_n}\nonumber\\
&\rightarrow& M
\end{eqnarray}
as $n\rightarrow +\infty$ for all $x\in(0,1)$.
It is obvious that (\ref{dnn}), (\ref{egn0}) and (\ref{egn1}) imply that $v_n:=\varphi_p\left(w_n'\right)$ is bounded in $C^1$. Therefore, by the Arzela-Ascoli theorem, we may assume that $v_n\rightarrow v$ in $C^0$. It follows that $w_n':=\varphi_q\left(v_n\right)\rightarrow \varphi_q(v)=:\widetilde{v}$ in $C^0$ where $q=p/(p-1)$. Obviously, $w_n(x)=\int_0^xw_n'(\tau)\,d\tau\rightarrow \int_0^x \widetilde{v}(\tau)\,d\tau$ in $C^0$.
Hence $w_n$ is strong convergence in $C^1$. Without loss of generality,
we may assume that $w_n\rightarrow w$ in $C^1$, $\Vert w\Vert=1$.
Clearly, we have $w\in \overline{S_k^\nu}$.

\emph{We claim that $w\in S_k^\nu$.}

On the contrary, suppose that $w\in \partial S_k^\nu$, then $w$ has at least one double zero $x_*\in[0,1]$. It follows that
$w_n\left(x_*\right)\rightarrow 0$ and $w_n'\left(x_*\right)\rightarrow 0$ as $n\rightarrow +\infty$. We note that
\begin{eqnarray}
\varphi_p\left(w_n'(x)\right)&=&\int_{x^*}^x\left(-\lambda_n a(\tau)\varphi_p\left(w_n\right)-\frac{f\left(\tau,u_n\left\vert u_n\right\vert^{\varepsilon_n},\lambda_n\right)}{\left\Vert u_n\right\Vert^{p-1}}-\frac{g\left(\tau,u_n,\lambda_n\right)}{\left\Vert u_n\right\Vert^{p-1}}\right)\,d\tau\nonumber\\
& &+\varphi_p\left(w_n'\left(x_*\right)\right).\nonumber
\end{eqnarray}
Firstly, we consider $x\in[0, x^*]$. Then we have that
\begin{eqnarray}
\left\vert w_n'(x)\right\vert&\leq&\varphi_p^{-1}\left(\left\vert\int_{x^*}^x\left(-\lambda_n a(\tau)\varphi_p\left(w_n\right)-\frac{F(\tau, u_n,\varepsilon_n,\lambda_n)}{\left\Vert u_n\right\Vert^{p-1}}\right)\,d\tau+\varphi_p\left(w_n'\left(x_*\right)\right)\right\vert\right)\nonumber\\
&\leq&\varphi_p^{-1}
\left(\int_x^{x^*}
 \left\vert\lambda_n a(\tau)\varphi_p\left(w_n\right)+\frac{F(\tau, u_n,\varepsilon_n,\lambda_n)}{\left\Vert u_n\right\Vert^{p-1}}\right\vert\,d\tau+\vert\varphi_p\left(w_n'\left(x_*\right)\right)\vert\right),\nonumber
\end{eqnarray}
where $F(\tau, u_n,\varepsilon_n,\lambda_n)=f\left(\tau,u_n\left\vert u_n\right\vert^{\varepsilon_n},\lambda_n\right)+g\left(\tau,u_n,\lambda_n\right)$.
It follows that
\begin{eqnarray}
\varphi_p(\left\vert w_n'(x)\right\vert)&\leq&\int_x^{x^*}
 \left\vert\lambda_n a(\tau)\varphi_p\left(w_n\right)+\frac{f\left(\tau,u_n\left\vert u_n\right\vert^{\varepsilon_n},\lambda_n\right)}{\left\Vert u_n\right\Vert^{p-1}}+\frac{g\left(\tau,u_n,\lambda_n\right)}{\left\Vert u_n\right\Vert^{p-1}}\right\vert\,d\tau\nonumber\\
 & &+\left\vert\varphi_p\left(w_n'\left(x_*\right)\right)\right\vert\nonumber\\
 &\leq&\int_x^{t^*}
 \left(\left\vert\lambda_n\right\vert a(\tau)+\left\vert\frac{f\left(\tau,u_n\left\vert u_n\right\vert^{\varepsilon_n},\lambda_n\right)}{\left\Vert u_n\right\Vert^{p-1}\varphi_p\left(w_n\right)}\right\vert+\left\vert\frac{g\left(\tau,u_n,\lambda_n\right)}{\left\Vert u_n\right\Vert^{p-1}\varphi_p\left(w_n\right)}\right\vert\right)\varphi_p\left(\left\vert w_n\right\vert\right)\,d\tau\nonumber\\
 & &+\left\vert\varphi_p\left(w_n'\left(x_*\right)\right)\right\vert.\nonumber
\end{eqnarray}
In view of (C2) and the definition of $w_n$, we can show that
\begin{equation}
\left\vert\frac{g\left(x,u_n,\lambda_n\right)}{\left\Vert u_n\right\Vert^{p-1}\varphi_p\left(w_n\right)}\right\vert=\left\vert\frac{g\left(x,u_n,\lambda_n\right)}{\varphi_p\left(u_n\right)}\right\vert\rightarrow 0\,\,\text{as}\,\,n\rightarrow+\infty\nonumber
\end{equation}
uniformly in $x\in(0,1)$.
Similarly, (C1) implies that
\begin{eqnarray}
\left\vert\frac{f\left(x,u_n\left\vert u_n\right\vert^{\varepsilon_n},\lambda_n\right)}{\left\Vert u_n\right\Vert^{p-1}\varphi_p\left(w_n\right)}\right\vert&=&\left\vert\frac{f\left(x,u_n\left\vert u_n\right\vert^{\varepsilon_n},\lambda_n\right)}{\varphi_p\left(u_n\left\vert u_n\right\vert^{\epsilon_n}\right)}\right\vert\varphi_p\left(\left\vert u_n\right\vert^{\epsilon_n}\right)\nonumber\\
&\leq& M\left\vert u_n\right\vert^{(p-1)\epsilon_n}\nonumber\\
&\rightarrow& M\nonumber
\end{eqnarray}
as $n\rightarrow +\infty$ for all $x\in(0,1)$. Hence there exists a positive constant $K$ such that
\begin{equation}
\left(\left\vert\lambda_n\right\vert a(x)+\left\vert\frac{f\left(x,u_n\left\vert u_n\right\vert^{\varepsilon_n},\lambda_n\right)}{\left\Vert u_n\right\Vert^{p-1}\varphi_p\left(w_n\right)}\right\vert+\left\vert\frac{g\left(x,u_n,\lambda_n\right)}{\left\Vert u_n\right\Vert^{p-1}\varphi_p\left(w_n\right)}\right\vert\right)\leq K\nonumber
\end{equation}
for all $x\in(0,1)$ and $n\in \mathbb{N}$ large enough. Thus, we have that
\begin{eqnarray}
\varphi_p\left(\left\vert w_n'(x)\right\vert\right)&\leq&K\int_x^{x^*}\varphi_p\left(\left\vert w_n\right\vert\right)\,d\tau+\vert\varphi_p\left(w_n'\left(x_*\right)\right)\vert.\nonumber
\end{eqnarray}
By the Gronwall-Bellman inequality [\ref{Bre}], we get that
$\varphi_p\left(\left\vert w_n'(x)\right\vert\right)\leq \left\vert\varphi_p\left(w_n'\left(x_*\right)\right)\right\vert \exp{K\left(x_*-x\right)}$. This means that $w_n'(x)\rightarrow 0$ on $[0, x^*]$ as $n\rightarrow+\infty$. Similarly, using a modification of Gronwall-Bellman inequality [\ref{ILL}, Lemma 2.2], we can get $w_n'(x)\rightarrow 0$ on $[x^*, 1]$ as $n\rightarrow+\infty$. It is obvious that $w_n(x)=\int_0^xw_n'(\tau)\,d\tau=w_n'(\xi)x\rightarrow 0$ as $n\rightarrow+\infty$, here $\xi\in[0,x]$. Therefore, $w_n\rightarrow 0$ in $C^1$ as $n\rightarrow+\infty$, which is a contradiction.

To obtain the bound on $\lambda$, we shall compare $w$ and $\psi_k^\nu$ via the generalized Picone identity [\ref{KJY}], where $\psi_k^\nu\in S_k^\nu$ is an eigenfunction of problem (\ref{pe}) corresponding to $\lambda_k$.
We have known that $w_n$ satisfies
\begin{equation}
\left(\varphi_p\left(w_n'\right)\right)'+\left(\lambda_n a(x)+\frac{f\left(x,u_n\left\vert u_n\right\vert^{\varepsilon_n},\lambda_n\right)}{\left\Vert u_n\right\Vert^{p-1}\varphi_p\left(w_n\right)}+\frac{g\left(x,u_n,\lambda_n\right)}{\left\Vert u_n\right\Vert^{p-1}\varphi_p\left(w_n\right)}\right)\varphi_p\left(w_n\right)=0\nonumber
\end{equation}
and $\psi_k^\nu$ satisfies
\begin{equation}
\left(\varphi_p\left(\left(\psi_k^\nu\right)'\right)\right)'+\lambda_k a(x)\varphi_p\left(\psi_k^\nu\right)=0.\nonumber
\end{equation}
Since $w_n$ and $\psi_k^\nu$ are both in $S_k^\nu$, by Lemma 2 of [\ref{Be}], there are two intervals $\left(\xi_1,\eta_1\right)$
and $\left(\xi_2,\eta_2\right)$ in $(0,1)$ where $w_n$ and $\psi_k^\nu$ do not vanish and have the same sign and such that $w_n\left(\xi_1\right)=w_n\left(\eta_1\right)=0$, and the same for $\left[\xi_2,\eta_2\right]$ with $w_n$ replaced by $\psi_k^\nu$.

We can assume without loss of generality that $w_n>0$ and $\psi_k^\nu>0$ in $\left(\xi_1,\eta_1\right)$. By the generalized Picone identity [\ref{KJY}],
we have that
\begin{eqnarray}\label{cp3}
-\int_{\xi_1}^{\eta_1}\left(\frac{w_n^{p-1}\varphi_p\left(\left(\psi_k^\nu\right)'\right)}{\varphi_p\left(\psi_k^\nu\right)}-w_n\varphi_p\left(w_n'\right)\right)'\,dx=A_1+B_1,
\end{eqnarray}
where
\begin{eqnarray}
A_1=\int_{\xi_1}^{\eta_1}\left(\lambda_ka(t)-\lambda_na(t)-\frac{f\left(t,u_n\left\vert u_n\right\vert^{\varepsilon_n},\lambda_n\right)}{\left\Vert u_n\right\Vert^{p-1}\varphi_p\left(w_n\right)}-\frac{g\left(t,u_n,\lambda_n\right)}{\left\Vert u_n\right\Vert^{p-1}\varphi_p\left(w_n\right)}\right)\varphi_p\left(w_n\right)\,dx\nonumber
\end{eqnarray}
and
\begin{eqnarray}
B_1=\int_{\xi_1}^{\eta_1}\left(\left\vert w_n'\right\vert^{p-1}+(p-1)\left\vert\frac{w_n\left(\psi_k^\nu\right)'}{\psi_k^\nu}\right\vert^p-p\varphi_p\left(w_n\right)w_n'\varphi_p\left(\frac{\left(\psi_k^\nu\right)'}{\psi_k^\nu}\right)\right)\,dx.\nonumber
\end{eqnarray}
The left-hand side of (\ref{cp3}) equals
\begin{equation}
\lim_{x\rightarrow \xi_1^+}\frac{w_n^p\varphi_p\left(\left(\psi_k^\nu\right)'\right)}{\varphi_p\left(\psi_k^\nu\right)}-
\lim_{x\rightarrow \eta_1^-}\frac{w_n^p\varphi_p\left(\left(\psi_k^\nu\right)'\right)}{\varphi_p\left(\psi_k^\nu\right)}:=H_{\xi_1}-H_{\eta_1}.\nonumber
\end{equation}
We prove that $H_{\xi_1}=0$. If $\psi_k^\nu\left(\xi_1\right)\neq 0$, then $H_{\xi_1}=0$. If $\psi_k^\nu\left(\xi_1\right)=0$, noting the conclusion of Lemma 2.2, then $\left(\psi_k^\nu\right)'\left(\xi_1\right)>0$.
By L'Hospital rule, we have that
\begin{eqnarray}
H_{\xi_1}&=&\lim_{x\rightarrow \xi_1^+}\frac{w_n^p\varphi_p\left(\left(\psi_k^\nu\right)'\right)}{\varphi_p\left(\psi_k^\nu\right)}\nonumber\\
&=&\lim_{x\rightarrow \xi_1^+}\frac{p\varphi_p\left(w_n\right)w_n'\varphi_p\left(\left(\psi_k^\nu\right)'\right)+w_n^p\left(\varphi_p\left(\left(\psi_k^\nu\right)'\right)\right)'}
{(p-1)\left(\psi_k^\nu\right)^{p-2}\left(\psi_k^\nu\right)'}\nonumber\\
&=&\lim_{x\rightarrow \xi_1^+}\frac{p\varphi_p\left(w_n\right)w_n'\varphi_p\left(\left(\psi_k^\nu\right)'\right)-w_n^p\lambda_k a(t)\varphi_p\left(\psi_k^\nu\right)}
{(p-1)\left(\psi_k^\nu\right)^{p-2}\left(\psi_k^\nu\right)'}.\nonumber
\end{eqnarray}
It follows that
\begin{eqnarray}
H_{\xi_1}&=&\lim_{x\rightarrow \xi_1^+}\frac{p\varphi_p\left(w_n\right)w_n'\varphi_p\left(\left(\psi_k^\nu\right)'\right)}
{(p-1)\left(\psi_k^\nu\right)^{p-2}\left(\psi_k^\nu\right)'}-\lim_{x\rightarrow \xi_1^+}\frac{w_n^p\lambda_k a(t)\varphi_p\left(\psi_k^\nu\right)}
{(p-1)\left(\psi_k^\nu\right)^{p-2}\left(\psi_k^\nu\right)'}\nonumber\\
&=&\lim_{x\rightarrow \xi_1^+}\frac{p\varphi_p\left(w_n\right)w_n'\varphi_p\left(\left(\psi_k^\nu\right)'\right)}
{(p-1)\left(\psi_k^\nu\right)^{p-2}\left(\psi_k^\nu\right)'}\nonumber\\
&=&\frac{pw_n'\left(\xi_1\right)\vert \left(\psi_k^\nu\right)'\left(\xi_1\right)\vert^{p-2}}{(p-1)}\lim_{x\rightarrow \xi_1^+}\frac{w_n^{p-1}}
{\left(\psi_k^\nu\right)^{p-2}}.\nonumber
\end{eqnarray}
If $p\leq2$, then $H_{\xi_1}=0$. If $2<p\leq 3$, applying L'Hospital rule again, we obtain that
\begin{equation}
\lim_{x\rightarrow \xi_1^+}\frac{w_n^{p-1}}{\left(\psi_k^\nu\right)^{p-2}}=\frac{(p-1)w_n'\left(\xi_1\right)}{(p-2)\left(\psi_k^\nu\right)'\left(\xi_1\right)}\lim_{x\rightarrow \xi_1^+}\frac{w_n^{p-2}}{\left(\psi_k^\nu\right)^{p-3}}.\nonumber
\end{equation}
This implies that $H_{\xi_1}=0$. If $k<p\leq k+1$, then we continue with this process $k$ times to obtain $H_{\xi_1}=0$.

Similarly, we can show that $H_{\eta_1}=0$. Therefore, the left-hand side of (\ref{cp3}) equals zero. Hence,
the right-hand side of (\ref{cp3}) also equals zero.

Young's inequality implies that
\begin{equation}
\left\vert w_n'\right\vert^{p-1}+(p-1)\left\vert\frac{w_n\left(\psi_k^\nu\right)'}{\psi_k^\nu}\right\vert^p-p\varphi_p\left(w_n\right)w_n'\varphi_p\left(\frac{\left(\psi_k^\nu\right)'}
{\psi_k^\nu}\right)\geq 0.\nonumber
\end{equation}
It follows that
\begin{eqnarray}\label{cp1}
A_1&=&\int_{\xi_1}^{\eta_1}\left(\lambda_ka(x)-\lambda_na(x)-\frac{f\left(x,u_n\left\vert u_n\right\vert^{\varepsilon_n},\lambda_n\right)}{\left\Vert u_n\right\Vert^{p-1}\varphi_p\left(w_n\right)}-\frac{g\left(x,u_n,\lambda_n\right)}{\left\Vert u_n\right\Vert^{p-1}\varphi_p\left(w_n\right)}\right)\varphi_p\left(w_n\right)\,dx\nonumber\\
&\leq& 0.
\end{eqnarray}
Similarly, we can also show that
\begin{eqnarray}\label{cp2}
A_2&=&\int_{\xi_2}^{\eta_2}\left(\lambda_na(x)-\lambda_ka(x)+\frac{f\left(x,u_n\left\vert u_n\right\vert^{\varepsilon_n},\lambda_n\right)}{\left\Vert u_n\right\Vert^{p-1}\varphi_p\left(w_n\right)}+\frac{g\left(x,u_n,\lambda_n\right)}{\left\Vert u_n\right\Vert^{p-1}\varphi_p\left(w_n\right)}\right)\varphi_p\left(\psi_k^\nu\right)\,dx\nonumber\\
&\leq& 0.
\end{eqnarray}

If $\lambda\leq \lambda_k$, considering (\ref{cp1}), (C1) and (C2), we have that
\begin{eqnarray}
\int_{\xi_1}^{\eta_1}\left(\lambda_ka(x)-\lambda a(x)\right)\varphi_p(w)\,dx&\leq& \lim_{n\rightarrow+\infty}\int_{\xi_1}^{\eta_1}\frac{f\left(x,u_n\left\vert u_n\right\vert^{\varepsilon_n},\lambda_n\right)}{\left\Vert u_n\right\Vert^{p-1}}\varphi_p\left(w_n\right)\,dx\nonumber\\
&\leq& \int_{\xi_1}^{\eta_1}M\varphi_p(w)\,dx.\nonumber
\end{eqnarray}
Hence, we get that
\begin{eqnarray}
\int_{\xi_1}^{\eta_1}\left(\lambda_k-\lambda \right)a_0\varphi_p(w)\,dx\leq \int_{\xi_1}^{\eta_1} M\varphi_p(w)\,dx,\nonumber
\end{eqnarray}
which implies $\lambda\geq \lambda_k-d$.

If $\lambda\geq\lambda_k$, considering (\ref{cp2}), (C1) and (C2), we have that
\begin{eqnarray}
\int_{\xi_2}^{\eta_2}\left(\lambda-\lambda_k\right)a(x)\varphi_p\left(\psi_k^\nu\right)\,dx&\leq& \lim_{n\rightarrow+\infty}\int_{\xi_2}^{\eta_2}\frac{-f\left(x,u_n\left\vert u_n\right\vert^{\varepsilon_n},\lambda_n\right)}{\left\Vert u_n\right\Vert^{p-1}\varphi_p\left(w_n\right)}\varphi_p\left(\psi_k^\nu\right)\,dx\nonumber\\
&\leq& \int_{\xi_2}^{\eta_2}M\varphi_p\left(\psi_k^\nu\right)\,dx.\nonumber
\end{eqnarray}
So, we obtain that
\begin{eqnarray}
\int_{\xi_2}^{\eta_2}\left(\lambda-\lambda_k\right)a_0\varphi_p\left(\psi_k^\nu\right)\,dx\leq \int_{\xi_2}^{\eta_2}M\varphi_p\left(\psi_k^\nu\right)\,dx,\nonumber
\end{eqnarray}
it follows $\lambda\leq\lambda_k+d$. Therefore, we have that $\lambda\in I_k$.
\qed
\\ \\
\textbf{Proof of Theorem 1.1.} We only prove the case of $\mathscr{D}_k^+$ since the proof of $\mathscr{D}_k^-$ can be given similarly.
Let $\mathscr{D}_k^+$ be the component of $\mathscr{S}_k^+\cup\left(I_k\times \{0\}\right)$ containing $I_k\times \{0\}$.
For any $(\lambda,u)\in \mathscr{D}_k^+$, there are two possibilities: (a) $u\in S_k^+$, or (b) $u\in \partial S_k^+$.
Clearly, $(\lambda,u)\in \Phi_k^+$ in the case of (a). While, the case (b) implies that $u$ has at least one double zero in $[0,1]$.
Lemma 2.2 follows that $u\equiv0$. Hence, there exists a sequence $\left(\lambda_n,u_n\right)\in \Phi_k^+$ such that
$\left(\lambda_n,u_n\right)$ is a solution of problem (\ref{dnf}) corresponding to $\epsilon=0$, and $\left(\lambda_n,u_n\right)$ converges to $(\lambda,0)$ in $\mathbb{R}\times E$.
By Lemma 2.3, we have $\lambda\in I_k$, i.e., $(\lambda,u)\in I_k\times \{0\}$ in the case of (b).
Hence, $\mathscr{D}_k^+\subset\Phi_k^+\cup\left(I_k\times \{0\}\right)$.

To complete the proof, it remains to show that $\mathscr{D}_k^+$ is unbounded in $\mathbb{R}\times E$.
Suppose on the contrary that $\mathscr{D}_k^+$ is bounded. Firstly, we claim that $\mathscr{D}_k^+$ is compact in $\mathbb{R}\times E$.

We consider the following auxiliary problem
\begin{equation}\label{eh}
\left\{
\begin{array}{l}
\left(\varphi_p\left(u'\right)\right)'=h,\,\,\text{a.e.}\,\, x\in (0,1),\\
u(0)=u(1)=0
\end{array}
\right.
\end{equation}
for a given $h\in L^1(0,1)$. It is known that problem (\ref{eh}) can be equivalently written as $u=G_p(h)(x)$. And $G_p:L^1(0,1)\to E$ is continuous and maps equi-integrable sets of $L^1(0,1)$ into relatively compacts of $E$. One may refer to
Lee and Sim [\ref{LS}] and Man\'{a}sevich and Mawhin [\ref{MM}] for
detail. Define the Nemitskii operator $H:\mathbb{R}\times E\rightarrow C(0,1)$
by
\begin{equation}
H(\lambda,u)(x):=-\lambda a(x)\varphi_p(u)-f\left(x,u,\lambda\right)-g(x,u,\lambda).\nonumber
\end{equation}
Then it is clear that $H$ is continuous operator which sends bounded
sets of $\mathbb{R} \times E$ into the equi-integrable sets of $C(0,1)$ and
problem (\ref{dF}) can be equivalently written as
\begin{equation}
u=G_p\circ H(\lambda,u):=\Psi_p(\lambda,u).\nonumber
\end{equation}
$\Psi_p$ is completely continuous in $\mathbb{R}\times E\rightarrow E$. Hence $\mathscr{D}_k^+$ is compact in $\mathbb{R}\times E$ because it
is bounded.

Applying a similar method to prove [\ref{Be}, Theorem 1] with obvious changes, we may find an open isolating neighborhood $\mathscr{O}$ of $\mathscr{D}_k^+$ such that $\partial \mathscr{O}\cap\mathscr{S}_k^+=\emptyset$. Note that such open isolating neighborhood can also be obtained by directly applying the results of [\ref{LM}, Proposition 3.2]
or [\ref{LM1}, Proposition 5.3] where the authors gave some general constructions of open isolating neighborhoods for Fredholm operators which can exhibit bifurcation from intervals.

In order to complete the proof of this theorem, we consider the approximate problem (\ref{dnf}) again. For $\epsilon>0$, it is easy to show that the nonlinear term $f\left(x,u\vert u\vert^\epsilon,\lambda\right)+g(x,u,\lambda)$ satisfies the condition (\ref{hg1}). Let
\begin{equation}
\mathscr{S}_\epsilon:=\overline{\left\{(\lambda,u)\big|(\lambda,u)\,\text{\,satisfies (\ref{dnf}) and $u\not\equiv 0$\,\,}\right\}}^{\mathbb{R}\times E}.\nonumber
\end{equation}
By Theorem 2.2 of [\ref{DM}], there exists a unbounded continuum $\mathscr{D}_{k,\epsilon}$ of $\mathscr{S}_\epsilon$ bifurcating from $\left(\lambda_k, 0\right)$, such that $\mathscr{D}_{k,\epsilon}\subset\Phi_k^\nu\cup\left\{\left(\lambda_k, 0\right)\right\}$. Furthermore, by Lemma 2.1, there are two sub-continua $\mathscr{D}_{k,\epsilon}^+$ and $\mathscr{D}_{k,\epsilon}^-$ of the continuum $\mathscr{D}_{k,\epsilon}$, which are both unbounded and $\mathscr{D}_{k,\epsilon}^\nu\subset\Phi_k^\nu\cup\left\{\left(\lambda_k, 0\right)\right\}$ or $\mathscr{D}_{k,\epsilon}^\nu\subset\Phi_k^{-\nu}\cup\left\{\left(\lambda_k, 0\right)\right\}$. Without loss of generality,
we may assume that
\begin{equation}
\mathscr{D}_{k,\epsilon}^\nu\subset\Phi_k^\nu\cup\left\{\left(\lambda_k, 0\right)\right\},\nonumber
\end{equation}
otherwise we can relabel it.

So there exists $\left(\lambda_\epsilon,u_\epsilon\right)\in \mathscr{D}_{k,\epsilon}^+\cap \partial \mathscr{O}$ for all $\epsilon>0$. By the compactness of $\Psi_p(\lambda,u)$, one can find a sequence $\epsilon_n\rightarrow 0$ such that
$\left(\lambda_{\epsilon_n},u_{\epsilon_n}\right)$ converges to a solution $(\lambda,u)$ of problem (\ref{dF}). So $u\in \overline{S_k^+}$.
If $u\in \partial S_k^+$, Lemma 2.2 or Remark 2.1 follows that $u\equiv 0$. By Lemma 2.3, $\lambda\in I_k$, which contradicts the definition of $\mathscr{O}$ (note that $(\lambda,u)\in \mathscr{D}_{k,\epsilon}^+\cap \partial \mathscr{O}$ since $\mathscr{D}_{k,\epsilon}^+\cap \partial \mathscr{O}$ is a closed subset of $\mathbb{R}\times E$). On the other hand, if $u\in S_k^+$, then $(\lambda,u)\in \mathscr{S}_k^+\cap \partial \mathscr{O}$ which contradicts $\partial \mathscr{O}\cap\mathscr{S}_k^+=\emptyset$. Therefore, $\mathscr{D}_k^+$ is unbounded in $\mathbb{R}\times E$.\qed\\

It is obvious that if $M=0$ the component $\mathscr{D}_k^+$ or $\mathscr{D}_k^-$ indeed
bifurcates from a single point $\lambda_k$. While, if $M\neq 0$, the following example shows that the component $\mathscr{D}_k^+$ or $\mathscr{D}_k^-$ possibly
bifurcates from an interval. \\
\\
\textbf{Example 2.1} (see [\ref{Be}]). Consider the following problem
\begin{equation}
\left\{
\begin{array}{l}
-u''=\lambda u+u\sin\left(u^2+u'^2\right)^{-1/2} \,\,\text{in}\,\,(0,\pi),\\
u(0)=u(\pi)=0.
\end{array}
\right.\nonumber
\end{equation}
It is easy to verify that $(\lambda(\rho),u(\rho)(x))=\left(1-\sin \vert \rho\vert^{-1},\rho\sin x\right)$ is a solution of this problem for any $\rho\neq 0$.
This follows that all the points of $[0,2]\times\{0\}$ are bifurcation points.
\\ \\
\textbf{Remark 2.2.} Note that two bifurcation intervals may be overlap. To see this, we consider the case of $a(x)\equiv1$ and $M>\left(2^p+1\right)\pi_p^p/2$,
where $\pi_p=\left(2\pi(p-1)^{1/p}\right)/
\left(p\sin\left(\pi/p\right)\right)$. In this case, it is well-known that $\lambda_k=\left(k\pi_p\right)^p$ (see [\ref{Z}]).
Then we can easily show that $\left(2\pi_p\right)^p-d<\pi_p^p+d$, which follows that $I_1\cap I_2\neq\emptyset$.
We can also see that two bifurcation intervals do not overlap if $k$ large enough. Anyway, $\mathscr{D}_k^\nu$ and $\mathscr{D}_j^\nu$ not meet
if $k\neq j$ since the nodal property.

\section{Spectrum of half-quasilinear problems}

\bigskip

\quad\, In this Section, we consider the half-quasilinear problem (\ref{hqle}).
Problem (\ref{hqle}) is called half-quasilinear because it is positive $p-$1-homogeneous and $p-$1-homogeneous in the cones $u>0$ and $u<0$.
Similar to that of [\ref{Be}], we say that $\lambda$ is a half-eigenvalue of problem (\ref{hqle}) if there exists a nontrivial solution $\left(\lambda,u_\lambda\right)$.
$\lambda$ is said to be simple if $v=cu_\lambda$, $c>0$ for all solutions ($\lambda,v$) of problem (\ref{hqle}).
\\ \\
\indent In order to prove Theorem 1.2, we need to establish Sturm type comparison theorem for $p$-Laplacian problems.\\
\\
\textbf{Lemma 3.1.} \emph{Let $b_2(x)\geq \max\left\{b_1(x),b_1(x)+\alpha-\beta\right\}$ for $x\in(0,1)$ and $b_i(x)\in C(0,1)$, $i=1,2$. Also let $u_1$, $u_2$
be solutions of the following differential equations
\begin{equation}
\left(\varphi_p\left(u'\right)\right)'+b_i(x)\varphi_p(u)+\alpha \varphi_p(u^+)+\beta \varphi_p(u^-)=0,\,\, i=1,2,\nonumber
\end{equation}
respectively. If $(c,d)\subset(0,1)$, and $u_1(c)=u_1(d)=0$, $u_1(x)\neq 0$ on $(c,d)$, then either there exists $\tau\in (c,d)$
such that $u_2(\tau)=0$ or $b_2=b_1$ and $u_2(x)=\mu u_1(x)$ for some constant $\mu\neq0$.}
\\ \\
\textbf{Proof.} If $u_2(x)\neq0$ on $(c,d)$, then we can assume without loss of generality that
$u_1(x)>0$, $u_2(x)>0$ in $(c,d)$. By the generalized Picone identity [\ref{KJY}] again, we have
\begin{eqnarray}\label{cp4}
-\int_c^d\left(\frac{u_1^p\varphi_p\left(u_2'\right)}{\varphi_p\left(u_2\right)}-u_1\varphi_p\left(u_1'\right)\right)'\,dx= \Gamma_1,
\end{eqnarray}
where
\begin{eqnarray}
\Gamma_1=\int_c^d\left(\left(b_2-b_1\right)u_1^p+\left(\left\vert u_1'\right\vert^p+(p-1)\left\vert \frac{u_1u_2'}{u_2}\right\vert^p-p\varphi_p\left(u_1\right)u_1'\varphi_p\left(\frac{u_2'}{u_2}\right)\right)\right)\,dx.\nonumber
\end{eqnarray}
The left-hand side of (\ref{cp4}) equals
\begin{equation}
\lim_{x\rightarrow c^+}\frac{u_1^p\varphi_p\left(u_2'\right)}{\varphi_p\left(u_2\right)}-u_1\varphi_p\left(u_1'\right)-
\lim_{x\rightarrow d^-}\frac{u_1^p\varphi_p\left(u_2'\right)}{\varphi_p\left(u_2\right)}-u_1\varphi_p\left(u_1'\right):=H_c-H_d.\nonumber
\end{equation}
By an argument similar to that of Lemma 2.3, we can show that $H_c=H_d=0$. Therefore, the left-hand side of (\ref{cp4}) equals to zero. Hence,
the right-hand side of (\ref{cp4}) also equals to zero.

Young's inequality implies that
\begin{equation}
\left\vert u_1'\right\vert^p+(p-1)\left\vert \frac{u_1u_2'}{u_2}\right\vert^p-p\varphi_p\left(u_1\right)u_1'\varphi_p\left(\frac{u_2'}{u_2}\right)\geq 0,\nonumber
\end{equation}
and the equality holds if and only if $\text{sgn}u_1'=\text{sgn} u_2'$ and $\left\vert u_1'/u_1\right\vert^p=\left\vert u_2'/u_2\right\vert^p$.
It follows that there exists a constant $\mu\neq0$ such that $u_2=\mu u_1$ and $b_2=b_1$.

Note that if $u_2(x)<0$ in $(c,d)$, similar to (\ref{cp4}), we can get
\begin{eqnarray}
-\int_c^d\left(\frac{u_1^p\varphi_p\left(u_2'\right)}{\varphi_p\left(u_2\right)}-u_1\varphi_p\left(u_1'\right)\right)'\,dx=\Gamma_2,\nonumber
\end{eqnarray}
where \begin{eqnarray}
\Gamma_2=\int_c^d\left(\left(b_2-b_1+\beta-\alpha\right)u_1^p+\left(\left\vert u_1'\right\vert^p+(p-1)\left\vert \frac{u_1u_2'}{u_2}\right\vert^p-p\varphi_p\left(u_1\right)u_1'\varphi_p\left(\frac{u_2'}{u_2}\right)\right)\right)\,dx.\nonumber
\end{eqnarray}
The above argument is still valid for this case.
\qed\\

By Lemma 3.1, we obtain the following result that will be used later.\\
\\
\textbf{Lemma 3.2.} \emph{Let $I_*=(a,b)$ be such $I_*\subseteq (0,1)$ and $\text{meas}\left\{I_*\right\}>0$.
Let $g_n:(0,1)\to \mathbb{R}$ be such that
\begin{equation}
\lim_{n\to +\infty} g_n(x)=+\infty\, \, \text{uniformly in}\,\,I_*.\nonumber
\end{equation}
Let $y_n\in E$ be a solution of the equation
\begin{equation}
\left(\varphi_p\left(y_n'\right)\right)'+ g_n(x)\varphi_p\left(y_n\right)+\alpha\varphi_p\left(y_n^+\right)+\beta\varphi_p\left(y_n^-\right)=0, \,\, x\in (0, 1).\nonumber
\end{equation}
 Then the number of zeros of $y_n\big|_{I_*}$ goes to infinity as $n\to +\infty$.
 }
\\ \\
\noindent{\bf Proof.} Set $\alpha^0:=\max_{x\in[0,1]}\vert\alpha(x)\vert$ and $\beta^0:=\max_{x\in[0,1]}\vert \beta(x)\vert$. By simple computation,
we can show that
\begin{equation}
g_n(x)+\alpha\varphi_p\left(\frac{y_n^+}{y_n}\right)+\beta\varphi_p\left(\frac{y_n^-}{y_n}\right)\geq g_n(x)-\alpha^0-\beta^0\,\,\text{for all}\,\, x\in(0,1).\nonumber
\end{equation}
After taking a subsequence if
necessary, we may assume that
\begin{equation}
g_{n_j}(x)-\alpha^0-\beta^0\geq j, \,\, x\in I_*\nonumber
\end{equation}
as $j\to +\infty$. It is easy to check the distance between any
two consecutive zeros of any nontrivial solution of the equation
\begin{equation}
\left\{
\begin{array}{l}
\left(\varphi_p\left(u'(x)\right)\right)'+j\varphi_p(u(x))=0,\,\, x\in I_*,\\
u(0)=u(1)=0
\end{array}
\right.\nonumber
\end{equation}
goes to zero as $j\to +\infty$. Note that the conclusion of Lemma 3.1 also is valid if $\alpha=\beta\equiv 0$. Using these facts
and Lemma 3.1, we can obtain the desired results.  \qed
\\ \\
\textbf{Proof of Theorem 1.2.} By Theorem 1.1, we know that there exists at least one solution of problem (\ref{hqle}), $\left(\lambda_k^\nu, u_k^\nu\right)\in \mathbb{R}\times S_k^\nu$,
for every $k=1,2,\ldots$, $\nu=+$ and $\nu=-$. The positive $p-$1-homogeneous of problem (\ref{hqle}) implies that
$\left\{\left(\lambda_k^\nu,cu_k^\nu\right),c>0\right\}$ are half-linear solutions in $\left\{\lambda_k^\nu\right\}\times S_k^\nu$.
Lemma 2.2 implies that any nontrivial solution $u$ of problem (\ref{hqle}) lies in some $S_k^\nu$.

\emph{We claim that for any solution $(\lambda,u)$ of problem (\ref{hqle}) with $u\in S_k^\nu$, we have that $\lambda=\lambda_k^\nu$ and $u=cu_k^\nu$ for some positive constant $c$.}

We may assume without loss of generality that the first zero of $u u_k^\nu$ to occur in $(0,1]$ is
a zero of $u$. That is, there exists $\zeta\in(0,1]$ such that
$u(\zeta)=0$, $u$ and $u_k^\nu$ do not vanish and have the same sign in $(0,\zeta)$. By Lemma
3.1 applied to $u$ and $u_k^\nu$ in $(0,\zeta)$, one has that $\lambda_k^\nu\leq \lambda$. On the other hand,
by Lemma 2 of [\ref{Be}], there must exist an interval $(\xi,\eta)\subset(0,1)$ such that $u$ and $u_k^\nu$
do not vanish and have the same sign in $(\xi,\eta)$, and $u_k^\nu(\xi)=u_k^\nu(\eta)=0$. Again by Lemma 3.1,
$\lambda\leq \lambda_k^\nu$; Hence $\lambda=\lambda_k^\nu$. Next, we shall prove that $u=cu_k^\nu$ for some positive constant $c$.

Without loss of generality, we may assume that $u$ and $u_k^\nu$ are positive in $(0,\zeta)$. Applying the generalized Picone identity to $u$ and $u_k^\nu$ on $[0,\zeta]$, noting $\lambda=\lambda_k^\nu$, we have that
\begin{eqnarray}\label{cp5}
-\int_0^\zeta\left(\frac{u^p\varphi_p\left(\left(u_k^\nu\right)'\right)}{\varphi_p\left(u_k^\nu\right)}-u\varphi_p\left(u'\right)\right)'\,dx=\Gamma_3,
\end{eqnarray}
where
\begin{eqnarray}
\Gamma_3=\int_0^\zeta\left(\left\vert u'\right\vert^p+(p-1)\left\vert \frac{u\left(u_k^\nu\right)'}{u_k^\nu}\right\vert^p-p\varphi_p(u)u'\varphi_p\left(\frac{\left(u_k^\nu\right)'}{u_k^\nu}\right)\right)\,dx.\nonumber
\end{eqnarray}
Using a proof similar to that of Lemma 2.3, we can show that the left-hand side of (\ref{cp5}) equals zero. Hence,
the right-hand side of (\ref{cp5}) also equals zero. It follows that $u=c_1u_k^\nu$ on $[0,\zeta]$ for some positive constant $c_1$.
We may assume without loss of generality that the first zero of $u u_k^\nu$ to occur in $(\zeta,1]$ is
a zero of $u$. That is, there exists $\zeta_1\in(\zeta,1]$ such that
$u(\zeta_1)=0$, $u$ and $u_k^\nu$ do not vanish and have the same sign in $(\zeta,\zeta_1)$.
Using method similar to the above, we can show that $u=c_2u_k^\nu$ on $[\zeta,\zeta_1]$ for some positive constant $c_2$.
Clearly, $u'(\zeta)=c_1\left(u_k^\nu\right)'(\zeta)=c_2\left(u_k^\nu\right)'(\zeta)$ and Lemma 2.2 imply $c_1=c_2$.
Repeating the above process $k$ times, we can show that $u=cu_k^\nu$ for some positive constant $c$.

By a similar method to prove [\ref{Be}, Theorem 2] with obvious changes, we can show that the sequences $\lambda_k^\nu$, $\nu=+$ or $-$ are increasing.\qed\\
\\
\textbf{Remark 3.1.} By simple computation, we can show that if $\beta\equiv 0$ then $\lambda_1^-=\lambda_1$, $\alpha\equiv 0$ implies
$\lambda_1^+=\lambda_1$ and $\alpha=\beta\equiv 0$ implies $\lambda_k^+=\lambda_k^-=\lambda_k$.
\\ \\
\textbf{Corollary 3.1.} \emph{If $\alpha=\beta$ in problem (\ref{hqle}), then $\lambda_k^+=\lambda_k^-:=\mu_k$ for each $k\in \mathbb{N}$.}
\\ \\
\textbf{Proof.} It is no difficulty to see that if the restriction of $c>0$ is replaced by $c\neq 0$, then the argument of Theorem 1.2 is also valid for
the case  $\alpha=\beta$. This fact combining with the conclusions of Theorem 1.2 implies the result. \qed
\\
\\
\textbf{Remark 3.2.} If $\alpha=\beta$, standard arguments by making use of the well-known Ljusternik-Schniremann theory [\ref{S}] and the techniques used in [\ref{A}],
we can obtain that problem (\ref{hqle}) has a sequence of eigenvalues
\begin{equation}
-\infty<\widetilde{\mu}_{1}<\widetilde{\mu}_{2}\leq\cdots\leq\widetilde{\mu}_{k}<\cdots,\lim_{k\rightarrow+\infty}\widetilde{\mu}_{k}=+\infty.\nonumber
\end{equation}
Corollary 3.1 implies that $\widetilde{\mu}_k=\mu_j$ for some $j\in \mathbb{N}$. Thus, $\widetilde{\mu}_k$ is simple.
And then Corollary 4.1 of [\ref{S}] implies that
\begin{equation}
-\infty<\widetilde{\mu}_{1}<\widetilde{\mu}_{2}<\cdots<\widetilde{\mu}_{k}<\cdots,\lim_{k\rightarrow+\infty}\widetilde{\mu}_{k}=+\infty.\nonumber
\end{equation}
Let $\sigma:=\left\{\mu_k\right\}_{k=1}^\infty$ and $\widetilde{\sigma}:=\left\{\widetilde{\mu}_k\right\}_{k=1}^\infty$. Clearly, $\widetilde{\sigma}\subseteq \sigma$.
We conjecture that $\widetilde{\sigma}=\sigma$. If $\alpha=\beta\leq \pi_p^p$ then by an argument similar to that of [\ref{ACM}, Theorem 1] with obvious changes,
we can show that $\widetilde{\sigma}=\sigma$. However, the methods used in [\ref{ACM}] become invalid for the general case $\alpha=\beta$.
\\

Naturally, we can consider the bifurcation structure of the perturbation of the problem (\ref{hqle}) of the form
\begin{equation}\label{hqle2}
\left\{
\begin{array}{l}
-\left(\varphi_p\left(u'\right)\right)'=\lambda a(x)\varphi_p(u)+\alpha \varphi_p\left(u^+\right)+\beta \varphi_p\left(u^-\right)+g(x,u,\lambda), \,\ x\in
(0,1),\\
u(0)=u(1)=0,
\end{array}
\right.
\end{equation}
where $g$ satisfies $\lim_{\vert u\vert\rightarrow 0}\left\vert g(x,u,\lambda)/\varphi_p(u)\right\vert=0$ uniformly in $x\in(0,1)$ and on bounded $\lambda$ intervals.\\
\\
\textbf{Theorem 3.1.} \emph{For each $k\in \mathbb{N}$, $\nu=+$ and $-$, $\left(\lambda_k^\nu,0\right)$ is a bifurcation point
for problem (\ref{hqle2}). Moreover, there exists a unbounded continuum $\mathscr{D}_k^\nu$ of solutions of problem (\ref{hqle2}), such that
$\mathscr{D}_k^\nu\subset\left(\left(\mathbb{R}\times S_k^\nu\right)\cup\left\{\left(\lambda_k^\nu,0\right)\right\}\right)$.}
\\ \\
\textbf{Proof.} Let us show that the only possible bifurcation points for problem (\ref{hqle2}) are the points $\left(\lambda_k^\nu,0\right)$.
Indeed, let $\left(\lambda_n, u_n\right)$, $u_n\not\equiv0$ be a sequence of solutions of problem (\ref{hqle2}) converging to $(\lambda,0)$. Let $v_n :=u_n/\left\Vert u_n\right\Vert$, then $v_n$ should be a solution of the problem
\begin{equation}
v_n= G_p\left( -\lambda_n a(x)\varphi_p\left(v_n(x)\right)-\alpha\varphi_p\left(v_n^+\right)-\beta\varphi_p\left(v_n^-\right)-\frac{g\left(x,u_n(x),\lambda\right)}{\left\Vert u_n(x)\right\Vert^{p-1}}\right).\nonumber
\end{equation}
By (\ref{egn0}) and the compactness of $G_p$, we obtain that for some convenient
subsequence $v_n\rightarrow v_0$ as $n\rightarrow+\infty$. Now $v_0$ verifies the equation
\begin{equation}
-\left(\varphi_p\left(v_0'\right)\right)'= \lambda a(x)\varphi_p(v_0)+\alpha\varphi_p(v_0^+)+\beta\varphi_p(v_0^-)\nonumber
\end{equation}
and $\left\Vert v_0\right\Vert = 1$. This implies that $\lambda=\lambda_k^\nu$ for some $k\in \mathbb{N}$ and $\nu\in\{+,-\}$.
The rest of proof is similar to that of [\ref{Be}, Theorem 3], so we omit it here.\qed\\
\\
\textbf{Remark 3.3.} Theorem 3.1 indicates that the bifurcation interval $I_k=\left\{\lambda_k^+,\lambda_k^-\right\}$, i.e.,
for problem (\ref{hqle2}), the bifurcation interval $I_k$ is a finite point set. What conditions can ensure that
the component indeed bifurcating from an interval is still an open problem even in the case of $p=2$.

\section{Nodal solutions for half-linear eigenvalue problems}

\bigskip

\quad\, We start this section by studying the following eigenvalue problem
\begin{equation}\label{hqle4}
\left\{
\begin{array}{l}
\left(\varphi_p\left(u'\right)\right)'+ \alpha \varphi_p\left(u^+\right)+ \beta \varphi_p\left(u^-\right)+\lambda r a(x)f(u)=0, \,\, x\in
(0,1),\\
u(0)=u(1)=0,
\end{array}
\right.
\end{equation}
where $\lambda>0$ is a parameter. Let $\zeta\in C(\mathbb{R})$ be such that
\begin{equation}
f(u)=f_0\varphi_p(u)+\zeta(u)\nonumber
\end{equation}
with $\lim_{\vert u\vert\rightarrow0}\zeta(u)/\varphi_p(u)=0$. Let us consider
\begin{equation}\label{hqle5}
\left\{
\begin{array}{l}
-\left(\varphi_p\left(u'\right)\right)'=\lambda r a(x)f_0\varphi_p(u)+ \alpha \varphi_p\left(u^+\right)+ \beta \varphi_p\left(u^-\right)+\lambda r a(x)\zeta(u),\,\, x\in
(0,1),\\
u(0)=u(1)=0
\end{array}
\right.
\end{equation}
as a bifurcation problem from the trivial solution $u\equiv 0$. \\

Applying Theorem 3.1 to problem (\ref{hqle5}), we have the following result.
\\ \\
\textbf{Lemma 4.1.} \emph{For each $k\in \mathbb{N}$, $\nu=+$ and $-$, $\left(\lambda_k^\nu/rf_0,0\right)$ is a bifurcation point of
problem (\ref{hqle5}). Moreover, there exists a unbounded continuum $\mathscr{D}_k^\nu$ of solutions of problem (\ref{hqle5}),
such that $\mathscr{D}_k^\nu\subset\Phi_k^\nu\cup\left\{\left(\lambda_k^\nu/rf_0,0\right)\right\}$.}
\\ \\
\textbf{Proof of Theorem 1.3.}
It is clear that any solution of problem (\ref{hqle4}) of the form $(1, u)$ yields a
solution $u$ of problem (\ref{dn}). We shall show that $\mathscr{D}_k^\nu$ crosses the hyperplane $\{1\}\times E$ in $\mathbb{R}\times E$.
To  this end, it will be enough to show that $\mathscr{D}_k^\nu$ joins
$\left(\lambda_k^\nu/rf_0, 0\right)$ to
$\left(\lambda_k^\nu/r f_\infty, +\infty\right)$. Let
$(\mu_n, y_n) \in \mathscr{D}_k^\nu\setminus\left\{\left(\lambda_k^\nu/rf_0,0\right)\right\}$ satisfy
\begin{equation}
\left\vert\mu_n\right\vert+\left\Vert y_n\right\Vert\rightarrow+\infty.\nonumber
\end{equation}

\emph{Case 1.}
$\lambda_k^\nu/f_\infty<r<\lambda_k^\nu/f_0$.

In this case, we only need to show that
\begin{equation}
\left(\frac{\lambda_k^\nu}{r f_\infty},\frac{\lambda_k^\nu}{r
f_0}\right) \subseteq\left\{\lambda\in\mathbb{R}\big|(\lambda,u)\in \mathscr{D}_k^\nu\right\}.\nonumber
\end{equation}
We divide the rest of proof into two steps.

\emph{Step 1.} We show that if there exists a constant $M>0$ such that $\left\vert\mu_n\right\vert\subset[0,M]$
for $n\in \mathbb{N}$ large enough, then $\mathscr{D}_k^\nu$ joins
$\left(\lambda_k^\nu/rf_0, 0\right)$ to
$\left(\lambda_k^\nu/r f_\infty, +\infty\right)$.

In this case it follows that
\begin{equation}
\left\Vert y_n\right\Vert\rightarrow+\infty.\nonumber
\end{equation}
Let $\xi\in C\left(\mathbb{R}\right)$ be such that
\begin{equation}
f\left(u\right)=f_\infty\varphi_p\left(u\right)+\xi\left(u\right)\nonumber
\end{equation}
with
\begin{equation}
\lim_{\vert u\vert\rightarrow+\infty}\frac{\xi\left(u\right)}{\varphi_p\left(u\right)}=0.\nonumber
\end{equation}
We divide the equation
\begin{equation}
-\left(\varphi_p\left(y_n'\right)\right)'=\mu_n r a(x)f_\infty \varphi_p\left(y_n\right)+ \alpha \varphi_p\left(y_n^+\right)+ \beta \varphi_p\left(y_n^-\right)+\mu_n r a(x)\xi\left(y_n\right)\nonumber
\end{equation}
by $\left\Vert y_n\right\Vert$ and set $\overline{y}_n =y_n/\left\Vert y_n\right\Vert$. Since $\overline{y}_n$ is bounded in $E$, after taking a subsequence if necessary, we have that $\overline{y}_n \rightharpoonup \overline{y}$ for some $\overline{y} \in E$ and $\overline{y}_n \rightarrow \overline{y}$ in $C[0,1]$.
Using the method to get (\ref{egn0}), we have that
\begin{equation}
\lim_{n\rightarrow+\infty}\frac{ \xi\left(y_n(x)\right)}{\left\Vert y_n\right\Vert^{p-1}}=0.\nonumber
\end{equation}
By the compactness of $G_p$, we obtain that
\begin{equation}
-\left(\varphi_p\left(\overline{y}'\right)\right)'=\overline{\mu} r a(x)f_\infty\varphi_p\left(\overline{y}\right)+ \alpha \varphi_p\left(\overline{y}^+\right)+ \beta \varphi_p\left(\overline{y}^-\right),\nonumber
\end{equation}
where $\overline{\mu}=\underset{n\rightarrow+\infty}\lim\mu_n$, again
choosing a subsequence and relabeling it if necessary.

It is clear that $\left\Vert \overline{y}\right\Vert=1$ and $\overline{y}\in
\overline{\mathscr{D}_{k}^\nu}\subseteq \mathscr{D}_{k}^\nu$ since
$\mathscr{D}_{k}^\nu$ is closed in $\mathbb{R}\times E$. Moreover, by Theorem 1.2, $\overline{\mu}r f_\infty=\lambda_k^\nu$, so
that
\begin{equation}
\overline{\mu}=\frac{\lambda_k^\nu}{r f_\infty}.\nonumber
\end{equation}
Therefore, $\mathscr{D}_{k}^\nu$ joins $\left(\lambda_k^\nu/r
f_0, 0\right)$ to $\left(\lambda_k^\nu/r f_\infty,
+\infty\right)$.

\emph{Step 2.} We show that there exists a constant $M$ such that $\left\vert\mu_n\right\vert
\in[0,M]$ for $n\in \mathbb{N}$ large enough.

On the contrary, we suppose that $\lim_{n\rightarrow +\infty}\left\vert\mu_n\right\vert=+\infty.$
Since $\left(\mu_n, y_n\right) \in \mathscr{D}_{k}^\nu$, it follows from the compactness of $G_p$ that
\begin{equation}
-\left(\varphi_p\left(y_n'\right)\right)'=\mu_n r \widetilde{f}_n(x)a(x)\varphi_p\left(y_n\right)+ \alpha \varphi_p\left(y_n^+\right)+ \beta \varphi_p\left(y_n^-\right),\nonumber
\end{equation}
where
\begin{equation}
\widetilde{f}_n(x)=\left\{
\begin{array}{l}
\frac{f\left(y_n(x)\right)}{\varphi_p\left(y_n(x)\right)},\,\, \text{if}\,\,y_n(x)\neq0,\\
f_0,\,\,\quad\quad\,\,\,\text{if}\,\,y_n(t)=0.
\end{array}
\right.\nonumber
\end{equation}
From (A1) and (A2), we can see that there exists a positive constant $\varrho$ such that
$f(s)/\varphi_p(s)\geq\varrho$ for all $s\neq 0$.
It follows that $\lim_{n\rightarrow+\infty}\mu_n r \widetilde{f}_n(x)=\pm\infty$. But if $\lim_{n\rightarrow+\infty}\mu_n r \widetilde{f}_n(x)=-\infty$, then $y_n\equiv 0$ on $(0,1)$ via Theorem 1.2, which is impossible. So $\lim_{n\rightarrow+\infty}\mu_n r \widetilde{f}_n(x)=+\infty$. By Lemma 3.2, we get that $y_n$ has more than $k-1$ zeros in $\left(0,1\right)$ for $n$
large enough, and this contradicts the fact that $y_n$ has exactly $k- 1$ zeros in $\left(0,1\right)$.

\emph{Case 2.} $\lambda_k^\nu/f_0<r<\lambda_k^\nu/f_\infty$.

Assume that $\left(\mu_n, y_n\right) \in \mathscr{D}_{k}^\nu$ is such that
\begin{equation}
\lim_{n\rightarrow+\infty}\left(\left\vert\mu_n\right\vert+\left\Vert y_n\right\Vert\right)=+\infty.\nonumber
\end{equation}
In view of \emph{Step 2} of \emph{Case 1}, we have known that there exists $M>0$, such that
for $n \in \mathbb{N}$ sufficiently large, $\left\vert\mu_n\right\vert\in [0,M].$
Applying the same method used in \emph{Step 1} of \emph{Case 1}, after
taking a subsequence and relabeling it if necessary, it follows that
\begin{equation}
\left(\mu_n,y_n\right)\rightarrow\left(\frac{\lambda_k^\nu}{r
f_\infty},+\infty\right) \,\, \text{as}\,\, n\rightarrow+\infty.\nonumber
\end{equation}
Thus, $\mathscr{D}_{k}^\nu$ joins $\left(\lambda_k^\nu/r
f_0,0\right)$ to $\left(\lambda_k^\nu/r
f_\infty,+\infty\right)$.\qed\\

\indent By an argument similar to that of Theorem 1.3, we can obtain the more general results as the following.
\\ \\
\textbf{Corollary 4.1.} \emph{Assume that $f$ satisfies (C1) and (A1), and for some $k, n\in\mathbb{N}$ with $k\leq n$, $\nu=+$ and $-$, either
$\lambda_n^\nu/f_{\infty}<r<\lambda_k^\nu/f_0$ or
$\lambda_n^\nu/f_0<r<\lambda_k^\nu/f_{\infty}.$ Then problem (\ref{dn}) has $n-k+1$ solutions $u_j^\nu$ for $j\in\{k,\cdots,n\}$ such that $u_j^\nu$
has exactly $j-1$ zeros in (0,1) and $\nu u_j^\nu$ is positive near 0.}
\\ \\
\indent Using Lemma 3.1, we can also get a nonexistence result when $f$ satisfies a non-resonant condition.
\\ \\
\textbf{Theorem 4.1.} \emph{Assume that $\alpha=\beta$ and there exists $k\in \mathbb{N}$ such that
\begin{equation}\label{ns}
\lambda_k^\nu<\frac{rf(s)}{\varphi_p(s)}<\lambda_{k+1}^\nu\,\, \text{for}\,\, s\neq 0 \text{\,\, and}\,\, r\in \mathbb{R}.
\end{equation}
Then problem (\ref{dn}) has no nontrivial solution in $\mathbb{R}\times E$.}
\\ \\
\textbf{Proof.} Suppose on the contrary that problem (\ref{dn}) has a nontrivial solution $u$ in $\mathbb{R}\times E$. Using the compactness of $G_p$, we can easily see that $u$ satisfies
\begin{equation}
\left(\varphi_p\left(u'\right)\right)'+\alpha \varphi_p\left(u^+\right)+\beta \varphi_p\left(u^-\right)+ rb(x)a(x)\varphi_p(u)=0,\,\,x\in (0,1),\nonumber
\end{equation}
where
\begin{equation}
b(x)=\left\{
\begin{array}{l}
\frac{f(u(x))}{\varphi_p(u(x))},\,\, \text{if}\,\,u(x)\neq0,\\
f_0,\,\,\,\quad\quad\text{if}\,\,u(x)=0.
\end{array}
\right.\nonumber
\end{equation}
If $u(x)\neq 0$, then (\ref{ns}) implies that $\lambda_k^\nu< rb(x)<\lambda_{k+1}^\nu$. On the other hand, for any zero $x_*$ of $u$,
we have $\lim_{x\rightarrow x_*}u(x)=0$. It follows that $\lim_{x\rightarrow t_*}f(u(x))/\varphi_p(u(x))=f_0$.
Hence, we have $b(x)\in C(0,1)$,
$$\lambda_k^\nu\leq rb(x)\leq\lambda_{k+1}^\nu,\,\,x\in (0,1)\nonumber$$
and all the inequalities are strict on some subset of positive measure in $(0,1)$.
We know that the eigenfunction $\varphi_k^\nu$ corresponding to $\lambda_k^\nu$ has exactly $k-1$ zeros in $(0,1)$.
Applying Lemma 3.1 to $\varphi_k^\nu$ and $u$, we can see that $u$ has at least $k$ zeros in $(0,1)$. By applying Lemma 3.1 again
to $u$ and $\varphi_{k+1}^\nu$, we can get that $\varphi_{k+1}^\nu$ has at least $k+1$ zeros in $(0,1)$. This is a contradiction.\qed\\
\\
\textbf{Remark 4.1.} Note that the assumption of $\alpha=\beta$ is crucial here to the application of Lemma 3.1. However, we doubt its necessity for this theorem.


\begin{thebibliography}{99}
\bibitem{}\label{A} A. Anane, Simplicit\'{e} et isolation de la premi\`{e}re valeur propre du $p$-Laplacien avec poids, C. R. Acad. Sci. Paris S\'{e}r.
I Math. 305 (1987), 725--728.

\bibitem{}\label{ACM} A. Anane, O. Chakrone and M. Monssa, Spectrum of one dimensional $p$-Laplacian with indefinite
weight, Electron. J. Qual. Theory Differ. Equ. 17 (2002), 1--11.

\bibitem{}\label{Be} H. Berestycki, On some nonlinear Sturm-Liouville
problems, J. Differential Equations 26 (1977), 375--390.

\bibitem{}\label{Bre} H. Brezis, Operateurs maximaux monotone et semigroup de contractions dans les espase de Hilbert,
Math. Studies, vol. 5, North-Holland, Amsterdam, 1973.

\bibitem{}\label{DM} G. Dai and R. Ma, Unilateral global bifurcation phenomena and nodal solutions
for $p$-Laplacian, J. Differential Equations 252 (2012), 2448--2468.

\bibitem{}\label{D1} E.N. Dancer, On the structure of solutions of non-linear eigenvalue problems, Indiana U. Math J.
23 (1974), 1069--1076.

\bibitem{}\label{D2} E.N. Dancer, Bifurcation from simple eigenvalues and eigenvalues of geometric multiplicity one, Bull.
Lond. Math. Soc. 34 (2002), 533--538.

\bibitem{}\label{DPM} M. Del Pino and R. Man\'{a}sevich, Global bifurcation from the
eigenvalues of the $p$-Laplacian, J. Differential Equations 92 (1991), 226--251.

\bibitem{}\label{DH} P. Dr\'{a}bek and Y,X. Huang, Bifurcation problems for the $p$-Laplacian in $\mathbb{R}^N$, Trans. Amer. Math. Soc. 349 (1997), 171--188.

\bibitem{}\label{E} L.C. Evans, Partial Differential Equations, AMS, Rhode Island, 1998.

\bibitem{}\label{FR} J. Fleckinger and W. Reichel, Global solution branches for $p$-Laplacian boundary value problems, Nonlinear Anal. 62 (2005), 53--70.

\bibitem{}\label{FMT} J. Fleckinger, R. Man\'{a}sevich and Th\'{e}lin, Global bifurcation from the first eigenvalue for
a system of $p$-Laplacians, Math. Nachr. 182 (1996), 217--242.

\bibitem{}\label{GS} J. Garc\'{\i}a-Meli\'{a}n and J. Sabina de Lis, A local bifurcation theorem for degenerate elliptic
equations with radial symmetry, J. Differential Equations 179 (2002), 27--43.

\bibitem{}\label{GT} P. Girg and P. Tak\'{a}\u{c}, Bifurcations of positive and negative continua
in quasilinear elliptic eigenvalue problems, Ann. Henri Poincar'e 9 (2008), 275--327.

\bibitem{}\label{H} J. K. Hale, Bifurcation from simple eigenvalues for several parameter families, Nonlinear Anal. 2 (1978), 491--497.

\bibitem{}\label{HK} P. Hess and T. Kato, On some linear and nonlinear eigenvalue problems with an indefinite weight
function, Comm. Partial Differential Equations 5 (1980), 999--1030.

\bibitem{}\label{ILL} B. Im, E. Lee and Y.H. Lee, A global bifurcation phenomena for second order singular boundary value problems,
J. Math. Anal. Appl. 308 (2005), 61--78.

\bibitem{}\label{K} M.A. Krasnosel'skii, Topological methods in the theory of nonlinear integral equations, Macmillan, New York, 1965.

\bibitem{}\label{KJY} T. Kusano, T. Jaros and N. Yoshida, A Picone-type identity and Sturmian comparison and oscillation theorems for a class of half-linear partial differential equations of second order, Nonlinear Anal. 40 (2000), 381--395.

\bibitem{}\label{LS} Y.H. Lee and I. Sim, Global bifurcation phenomena for singular
one-dimensional $p$-Laplacian, J. Differential Equations 229 (2006), 229--256.

\bibitem{}\label{L0} J. L\'{o}pez-G\'{o}mez, Multiparameter bifurcation based on the linear part, J. Math. Anal. Appns. 138 (1989),
358--370.

\bibitem{}\label{L2} J. L\'{o}pez-G\'{o}mez, Positive periodic solutions of Lotka-Volterra RD Systems, Diff. Int. Eqns. 5 (1992),
55--72.

\bibitem{}\label{L1} J. L\'{o}pez-G\'{o}mez, Spectral theory and nonlinear functional analysis, Chapman and Hall/CRC, Boca
Raton, 2001.

\bibitem{}\label{LM} J. L\'{o}pez-G\'{o}mez and C. Mora-Corral, Minimal complexity of semi-bounded components in bifurcation
theory, Nonlinear Anal. 58 (2004), 749--777.

\bibitem{}\label{LM1} J. L\'{o}pez-G\'{o}mez and C. Mora-Corral, Counting zeroes of $C^1$-Fredholm maps of index 1, Bull. Lond.
Math. Soc. 37 (2005), 778--792.

\bibitem{}\label{LM2} J. L\'{o}pez-G\'{o}mez and C. Mora-Corral, Algebraic Multiplicity of Eigenvalues of Linear Operators,
Advances in Operator Theory and Applications Vol. 177, Birkha\"{u}ser, Basel, 2007.

\bibitem{}\label{MT1} R. Ma and B. Thompson, Nodal solutions for nonlinear eigenvalue problems,
Nonlinear Anal. 59 (2004), 707--718.

\bibitem{}\label{MM} R. Man\'{a}sevich and J. Mawhin,  Periodic solutions for nonlinear
systems with $p$-Laplacian-like operators, J. Differential Equations
145 (1998), 367--393.

\bibitem{}\label{R1} P.H. Rabinowitz, Nonlinear Sturm-Liouville problems for second order ordinary
differential equations, Commun. Pure Appl. Math. 23 (1970), 939--961.

\bibitem{}\label{R2} P.H. Rabinowitz, Some global results for nonlinear eigenvalue problems,
J. Funct. Anal. 7 (1971), 487--513.

\bibitem{}\label{R3} P.H. Rabinowitz, On bifurcation from infinity,
J. Funct. Anal. 14 (1973), 462--475.

\bibitem{}\label{R4} P.H. Rabinowitz, Some aspects of nonlinear eigenvalue problems,
Rocky Mountain J. Math. 3 (1973), 161--202.

\bibitem{}\label{SS} K. Schmitt and H.L. Smith, On eigenvalue problems for nondifferentiable mappings,
J. Differential Equations 33 (1979), 294--319.

\bibitem{}\label{S} A. Szulkin, Ljusternik-Schnirelmann theory on $C^1$-manifolds, Ann. Inst. H. Poincar\'{e} Anal. Non Lin\'{e}aire 5 (1988), 119--139.

\bibitem{}\label{Z} M.R. Zhang, The rotation number approach to eigenvalues of the one-dimensional
$p$-Laplacian with periodic potentials, J. Lond. Math. soc. (2) 64 (2001), 125--143.
\end{thebibliography}
\end{document}